# On a new class of fractional partial differential equations II

## Tien-Tsan Shieh and Daniel E. Spector

*Dedicated to Irene Fonseca, on the occasion of her 60th birthday, with esteem and affection*

**Abstract.** In this paper we continue to advance the theory regarding the Riesz fractional gradient in the calculus of variations and fractional partial differential equations begun in an earlier work of the same name. In particular we here establish an $L^1$ Hardy inequality, obtain further regularity results for solutions of certain fractional PDE, demonstrate the existence of minimizers for integral functionals of the fractional gradient with non-linear dependence in the field, and also establish the existence of solutions to corresponding Euler-Lagrange equations obtained as conditions of minimality. In addition we pose a number of open problems, the answers to which would fill in some gaps in the theory as well as to establish connections with more classical areas of study, including interpolation and the theory of Dirichlet forms.

## 1 Introduction

In the preceding paper of the same name [42], the authors undertook the exposition of the Riesz fractional gradient and its systematic study from the perspective of the calculus of variations. Here we recall that for $s \in (0, 1)$ one can define in $d$-dimensional Euclidean space the fractional gradient by

$$D^s u(x) := c_{d,s} \int_{\mathbb{R}^d} \frac{u(x) - u(y)}{|x-y|^{d+s}} \frac{x-y}{|x-y|} \, dy, \tag{1}$$

The first author is partially supported by National Science Council of Taiwan under research grant NSC 101-2115-M-009-014-MY3. The second author is supported by the Taiwan Ministry of Science and Technology under research grants 103-2115-M-009-016-MY2 and 105-2115-M-009-004-MY2.

provided $u$ has sufficient smoothness and decay such that (1) is well-defined as a Lebesgue integral. The fractional gradient is the canonical example[1] of the non-local gradients considered by Mengesha and the second author in [33], where localization results were obtained for these non-local gradients and integral functionals defined in terms of them, while it can be contrasted with the more well-known fractional Laplacian

$$(-\Delta)^{s/2} u(x) := \tilde{c}_{d,s} \int_{\mathbb{R}^d} \frac{u(x) - u(y)}{|x - y|^{d+s}} \, dy \qquad (2)$$

as a curl free vector with the same differential order. The latter has been the subject of extensive study, making its way into the canon of literature in both harmonic analysis [45] and the study of fractional derivatives [34], while recently there has been a renewed interest in it from the standpoint of fractional partial differential equations [5, 9–11, 14, 15, 28, 30, 31, 35].

As the fractional gradient has not received such prominent attention, the purpose of the preceding paper was to introduce it as a fundamental object of study. In particular, we showed that with such a definition one can continuously interpolate the class of minimization problems in the calculus of variations with linear dependence in the field from differential order zero to one: For $\Omega \subset \mathbb{R}^d$ open, we established existence of minimizers of integral functionals of the form

$$F_s(u) = \int_{\mathbb{R}^d} f(x, D^s u) - g u \chi_\Omega \, dx, \qquad (3)$$

for suitable hypothesis on $f, g$. We further showed that such minimizers are solutions to corresponding Euler-Lagrange equations, a form of fractional partial differential equations arising from the conditions of minimality.

**Remark 1.1.** A consensus as to what constitutes a general fractional elliptic equation has not been made, though any candidate should contain the fractional Laplace's equation as its most basic example. Various theories to this effect have been pursued in a number of directions, see for example a "fully nonlinear divergence form theory" of Caffarelli and Silvestre [9, 10], a "divergence form elliptic complex interpolation theory" of Caffarelli and Stinga [11], and a "divergence form elliptic real interpolation theory" with contributions by a number of authors: Brasco and Lindgren [5], Di Castro, Kuusi, and Palatucci [14, 15], Korvenpaa, Kuusi, and Palatucci [28], Kuusi, Mingione and Sire [30, 31] and Schikorra [39].

---

[1]   Technically the non-local gradients in [33] were defined with integration over a bounded domain, though with suitable modification much of the analysis can be extended to integration over all of space.

As far as the authors are aware, there was no formal name for the object (1) preceding our paper, though it may be recognizable[2] through the relation

$$D^s u(x) \equiv R(-\Delta)^{s/2} u(x) \qquad (4)$$

for $R$ the vector valued Riesz transform:

$$Rf(x) = c_d \ p.v. \int_{\mathbb{R}^d} f(y) \frac{x-y}{|x-y|^{d+1}} \ dy.$$

In this respect it is prudent to make here a more thorough review of papers utilizing the fractional gradient that have come to our attention since the publication of [42]. The earliest reference we are aware of concerning an equivalent definition to (1) is the 1959 paper of Horváth [23], while it is implicit in the 1961 paper of Sobolev and Nikol'skiĭ- see p. 148 in [43]. From the standpoint of applications, a bounded domain analogue of (1) can be subsumed in the non-local continuum mechanics theory developed by Edelen and Laws [17, Equation (3.1), p. 27] and Edelen, Green and Laws [16, Equation (3.1), p. 38]. In these several papers, the standard local PDE - local conservation of mass, balance of momentum, and balance of moment of momentum- are equipped with a global balance of energy and global variational postulate in the constitutive equations and global balance of energy and global Clausius-Duhem inequality, respectively. In particular, the constitutive assumption in Equation (3.2) in [17] assumes the internal energy depends on the non-local substate variables which in turn depend on the deformation gradient from the deduction following Equation (3.8) on p. 28, while in Equation (3.2) in [16] the Helmholtz free energy is assumed to depend on a similar non-local substate variable depending on the deformation gradient. By taking this non-local substate variable to be a convolution with a restricted Riesz potential, one finds a local version of (1). In more contemporary work it has appeared in the papers of Caffarelli and Vazquez [7], Caffarelli, Soria, and Vazquez [6] and Biler, Imbert and Karch [3] in the context of a non-local porous medium equation as the gradient of the Riesz potential of the density - a "fractional potential pressure" (in particular the formula (1) has been recorded in [3]). The fractional gradient's appearance as a boundary-type operator in the spirit of Caffarelli and Silvestre's result [8] has been obtained by Stinga and Torrea in [46], while in [36, 37] Schikorra has considered a vector-valued analogue in the context of harmonic maps, establishing regularity for critical points of conformal energies of the fractional gradient. Finally let us mention a second order analogue - a fractional Hessian - has been considered by Guillen and Schwab in [22].

---

[2] L. Caffarelli gave the second author such a description in a discussion in 2014 in Haifa, Israel.

As developed in the preceding paper, a motivation for considering either (1) or (2) can be found in the desire for a theory which gives rise to spaces with good functional properties - compactness, embeddings, etc. Indeed this was precisely the aim of Sobolev and Nikol'skiĭ's paper [43]. In fact, in one dimension both (1) and (2) can be obtained explicitly from linear combinations of the Liouville fractional derivatives they suggest for such a theory:

$$\frac{d^s}{dx_+^s} u(x) = \frac{1}{\Gamma(-s)} \int_0^\infty \frac{u(x+h) - u(x)}{|h|^{1+s}} \frac{h}{|h|} \, dh$$

and

$$\frac{d^s}{dx_-^s} u(x) = \frac{-1}{\Gamma(-s)} \int_{-\infty}^0 \frac{u(x+h) - u(x)}{|h|^{1+s}} \frac{h}{|h|} \, dh.$$

In particular,

$$(\frac{d^s}{dx_+^s} - \frac{d^s}{dx_-^s})u(x) = cD^s u(x)$$

$$(\frac{d^s}{dx_+^s} + \frac{d^s}{dx_-^s})u(x) = \tilde{c}(-\Delta)^{s/2}u(x),$$

while the relation (4) also holds in this one dimensional example with the Hilbert transform in place of the Riesz transforms.

The purpose of this paper is to continue to advance the theory begun in the paper [42] and built upon in the subsequent papers developing inequalities for the fractional gradient [38] and regularity of solutions of fractional PDE defined in terms of it [40, 41]. One aspect of this development is to establish a number of new results for the fractional gradient. This includes an $L^1$ Hardy inequality, further regularity results for fractional PDE, the existence of minimizers for certain integral functionals of the fractional gradient with nonlinear dependence in the field, and also the existence of solutions to corresponding Euler-Lagrange equations obtained as conditions of minimality. Yet despite the several papers we have written and the current advances in the theory we present here, given the general lack of study of the fractional gradient as a fundamental object, there is still much to be explored. Therefore a second aspect of this paper is to present some open problems we have formulated in the course of our research. The answers to some of these questions would fill in details currently missing in our understanding of fractional phenomena that have been settled in the integer setting. The answer to others are of interest because they would establish connections with areas of classical interest such as complex interpolation or the theory of Dirichlet forms. In any case, answers to these questions would certainly provide new insight and tools that would be useful in future work.

## 1.1 An $L^1$ Hardy Inequality for the Fractional Gradient

The first item we address in this paper is the existence of an $L^1$ Hardy inequality for the fractional gradient. This question of inequalities for the fractional gradient in the $L^1$ regime was taken up in [38] in the case of Sobolev's inequality. Let us here recall that while in this endpoint the Sobolev inequality fails for the fractional Laplacian (cf. [45, p.119]), in [38] the following result was demonstrated: For any $s \in (0,1)$ there exists a constant $C = C(s,d) > 0$ such that

$$\|u\|_{L^{d/(d-s)}(\mathbb{R}^d)} \leq C \|D^s u\|_{L^1(\mathbb{R}^d; \mathbb{R}^d)},$$

for all $u$ such that $D^s u \in L^1(\mathbb{R}^d; \mathbb{R}^d)$. One observes a similar difficulty for a fractional Hardy inequality in the $L^1$ endpoint, since while it is a simple consequence of a result of Stein and Weiss [44] that for $1 < p < +\infty$ one has the inequality

$$\int_{\mathbb{R}^d} \frac{|u|^p}{|x|^{sp}} \, dx \leq C \int_{\mathbb{R}^d} |D^s u|^p \, dx, \tag{5}$$

the same counterexample as in Sobolev's inequality shows that one cannot have (5) when $p = 1$ with the fractional Laplacian on the right hand side. Nonetheless, the first result we show in this paper is

**Theorem 1.2.** *For all $s \in (0,1)$ one has*

$$\frac{(d-1)\Gamma(\frac{s}{2})\Gamma(\frac{d-1}{2})}{\pi^{(2-s)/2}2^{1-s}\Gamma(\frac{d-s}{2})} \int_{\mathbb{R}^d} \frac{|u|}{|x|^s} \, dx \leq \int_{\mathbb{R}^d} |D^s|u|| \, dx.$$

*for all $u$ such that $D^s|u| \in L^1(\mathbb{R}^d; \mathbb{R}^d)$.*

The argument is surprisingly simple, as it follows the proof of the classical Hardy inequality. In fact, we would have such an inequality for all $u$ such that $D^s u \in L^1(\mathbb{R}^d; \mathbb{R}^d)$ if one could answer

**Open Problem 1.3.** Does there exist a constant $C > 0$ (possibly depending on $d$ and $s$) such that

$$\int_{\mathbb{R}^d} |D^s|u|| \, dx \leq C \int_{\mathbb{R}^d} |D^s u| \, dx$$

for all $u$ such that $D^s u \in L^1(\mathbb{R}^d; \mathbb{R}^d)$?

The answer to such a question for $1 < p < +\infty$ follows easily from complex interpolation, since one has (see, for example, Bergh and Löfström [2, p. 153, (7)])

$$(L^p(\mathbb{R}^d), W^{1,p}(\mathbb{R}^d))_{[s]} = H^{s,p}(\mathbb{R}^d),$$

which combined with the sublinearity of the map $z \mapsto |z|$ and its boundedness on $L^p(\mathbb{R}^d)$ and $W^{1,p}(\mathbb{R}^d)$ allows one to invoke the result of Calderón and Zygmund [12]. When $p = 1$ the reliance of the above interpolation argument on retracts would yield interpolation of the Hardy space $\mathcal{H}^1(\mathbb{R}^d)$ and a Hardy-Sobolev space, the result of which would be an estimate involving both the fractional gradient and the fractional Laplacian. Yet such an inequality can already be deduced from the $L^1$ Hardy inequality utilizing Gagliardo semi-norms (see, e.g. Frank and Seiringer [19]) as the space $W^{s,1}(\mathbb{R}^d)$ can be seen to be embedded in the space of all locally integrable functions whose fractional Laplacian is in $\mathcal{H}^1(\mathbb{R}^d)$ by integrating the basic inequality

$$|D^s u(x)| + |(-\Delta)^{s/2} u(x)| \le C \int_{\mathbb{R}^d} \frac{|u(x) - u(y)|}{|x - y|^{d+s}} \, dy,$$

which follows from the definitions (1) and (2). Given this difficulty when $p = 1$, one is interested to understand if the space of functions such that $D^s u \in L^1(\mathbb{R}^d; \mathbb{R}^d)$ is an interpolation space. In fact, one wonders

**Open Problem 1.4.** Is $W^{1,1}(\mathbb{R}^d)$ an interpolation space?

Even considering the real interpolation of $L^1(\mathbb{R}^d)$ and $\dot{W}^{2,1}(\mathbb{R}^d)$ one has (cf. [2, p. 147])

$$\left( L^1(\mathbb{R}^d), \dot{W}^{2,1}(\mathbb{R}^d) \right)_{1/2,q} = B^1_{1q}(\mathbb{R}^d)$$

where $B^1_{1q}(\mathbb{R}^d)$ are Besov spaces and not the classical Sobolev space $W^{1,1}(\mathbb{R}^d)$.

## 1.2 Regularity for fractional PDE

Secondly, as was observed in the recent papers of the authors in collaboration with Armin Schikorra [40, 41], the question of regularity for fractional partial differential equations in this framework follows the classical regularity. For example, in [40] we extended the technique of Iwaniec and Sbordone [24] to obtain the following result: Let $\Omega \subset \mathbb{R}^d$ be open. Suppose $A : \mathbb{R}^d \to \mathbb{R}^{d \times d}$ is a function of vanishing mean oscillation and uniformly elliptic, i.e.

$$\sup_Q \fint_Q \left| A - \fint_Q A \right| \, dx < +\infty,$$

$$\lim_{|Q| \to 0, \infty} \fint_Q \left| A - \fint_Q A \right| \, dx = 0$$

and

$$\lambda |\xi|^2 \leq A(x)\xi \cdot \xi \leq \Lambda |\xi|^2, \tag{6}$$

for all $x, \xi \in \mathbb{R}^d$ and some $0 < \lambda \leq \Lambda < +\infty$. Further suppose $G \in L^p(\mathbb{R}^d; \mathbb{R}^d)$ and

$$\int_{\mathbb{R}^d} A(x) D^s u(x) \cdot D^s \varphi(x) \, dx = \int_{\mathbb{R}^d} G \cdot D^s \varphi \tag{7}$$

for all $\varphi \in C_c^\infty(\Omega)$. Then $D^s u \in L_{loc}^p(\Omega)$ and for any $K \subset\subset \Omega$ there exists a constant $C = C(K, \Omega, A, s, p) > 0$ such that

$$\|D^s u\|_{L^p(K; \mathbb{R}^d)} \leq C \left( \|G\|_{L^p(\mathbb{R}^d; \mathbb{R}^d)} + \|(-\Delta)^{s/2} u\|_{L^2(\mathbb{R}^d)} \right).$$

More recently, in [41] we introduced a reduction argument that amounts to lifting the fractional PDE to a classical equation to obtain regularity of the homogenous equation for an $H^{s,p}$-Laplacian: Let $\Omega \subset \mathbb{R}^d$ be open, $p \in (2 - \frac{1}{d}, \infty)$ and $s \in (0, 1]$. Suppose $u \in H^{s,p}(\mathbb{R}^d)$ satisfies

$$\int_{\mathbb{R}^d} |D^s u|^{p-2} D^s u \cdot D^s \varphi = 0 \quad \forall \varphi \in C_c^\infty(\Omega). \tag{8}$$

Then $u \in C_{loc}^{s+\alpha}(\Omega)$ for some $\alpha > 0$ only depending on $p$.

In fact, this technique extends to a larger class of homogeneous equations for which regularity is known in the classical inhomogeneous case, a result we now develop. Following [29], we assume that $a : \mathbb{R}^d \times \mathbb{R}^d \to \mathbb{R}^d$ satisfies the growth, ellipticity, and continuity assumptions:

$$|a(x, \xi)| + |\partial_\xi a(x, \xi)|(|\xi|^2 + s^2)^{1/2} \leq L(|\xi|^2 + s^2)^{(p-1)/2}$$
$$\nu(|\xi|^2 + s^2)^{(p-2)/2}|\zeta|^2 \leq \partial_\xi a(x, \xi)\zeta \cdot \zeta \tag{9}$$

for all $x, \xi, \zeta \in \mathbb{R}^d$; $\partial_\xi a$ is assumed to be continuous in $\xi$ if $p \geq 2$ and continuous away from the origin if $p < 2$; $a$ is assumed to be measurable in $x$; $\nu, L, s$ are fixed parameters with $0 < \nu \leq L$ and $s \geq 0$. Then a function $v$ is a weak solution to the equation

$$- \operatorname{div} a(x, Dv) = \mu$$

if $v \in W^{1,p}_{loc}(\Omega)$ and

$$\int_{\Omega} a(x, Dv) \cdot D\varphi = \int_{\Omega} \varphi \, d\mu \tag{10}$$

for all $\varphi \in C_c^{\infty}(\Omega)$. Further define the averaged renormalized modulus of continuity of $x \mapsto a(x, \cdot)$ as

$$\omega(r) := \left[ \sup_{z \in \mathbb{R}^d, B(x,r) \subset \Omega} \fint_{B(x,r)} \left( \frac{|a(y, \xi) - \fint_{B(x,r)} a(w, \xi) \, dw|}{(|\xi|^2 + s)^{p-1}} \right)^2 \, dy \right]^{1/2}.$$

If $p \geq 2$, we assume the Dini-Hölder condition

$$S := \sup_r \int_0^r \frac{[\omega(\rho)]^{2/p}}{\rho^{\tilde{\alpha}}} \, \frac{d\rho}{\rho} < +\infty \tag{11}$$

for some $\tilde{\alpha} < \alpha_M$, the maximal Hölder regularity of $Dv$ satisfying the homogeneous equation

$$- \operatorname{div} a(Dv) = 0.$$

If $p \in (2 - \frac{1}{d}, 2]$ we assume that

$$S := \sup_r \int_0^r \frac{[\omega(\rho)]^{\sigma}}{\rho^{\tilde{\alpha}}} \, \frac{d\rho}{\rho} < +\infty \tag{12}$$

for some $\sigma < 1$. Then by Theorems 1.4 and 1.6 in [29], we find that for $\mu \in L^\infty_{loc}(\Omega)$ and $v$ satisfying (10) one has $Dv \in C^\alpha_{loc}(\Omega)$. We here apply the same reduction argument as in [41] to obtain the following regularity result.

**Theorem 1.5.** *Suppose that $p \in (2 - \frac{1}{d}, \infty)$ and that $a(x, \xi)$ satisfies (9). If $p \geq 2$, further assume that $a$ satisfies (11), while if $p < 2$ assume that $a$ satisfes (12). Moreover, we additionally assume that for all $|x|$ sufficiently large $|a(x, \xi)| \leq L|\xi|^{p-1}$. Then for any $u \in H^{s,p}_g(\Omega)$ that satisfies*

$$\int_{\mathbb{R}^d} a(x, D^s u) \cdot D^s \varphi \, dx = 0$$

*for all $\varphi \in C_c^{\infty}(\Omega)$, one has $u \in C^{s+\alpha}_{loc}(\Omega)$.*

Returning to the linear equation when $A$ is only assumed to be bounded, measurable, and satisfy (6), this argument allows us to show

**Theorem 1.6.** *Suppose $A : \mathbb{R}^d \to \mathbb{R}^{d \times d}$ is bounded, measurable, and satisfies* (6). *Further suppose $G \in L^p(\mathbb{R}^d; \mathbb{R}^d)$ for some $p > d/s$ and*

$$\int_{\mathbb{R}^d} A(x) D^s u \cdot D^s \varphi(x) \, dx = \int_{\mathbb{R}^d} G \cdot D^s \varphi \, dx.$$

*for all $\varphi \in C_c^\infty(\Omega)$. Then $I_{1-s} u \in C_{loc}^\alpha(\Omega)$.*

However, this regularity does not match that obtained in [41] for $A \in VMO$, and so one wonders

**Open Problem 1.7.** Is it true that for every $s \in (0, 1)$, any $u$ satisfying (7) with $G \in L^p(\mathbb{R}^d; \mathbb{R}^d)$ for some $p > d/s$ and $A$ only bounded and measurable is Hölder continuous of some exponent $\alpha > 0$ (which possibly depends upon $s$)?

This can be compared with the classical setting, where it was De Giorgi [13] who proved the Hölder regularity of solutions to elliptic equations with bounded and measurable coefficients. Let us recall his result here, which for convenience of reference and exposition we follow the formulation of Kinderlehrer and Stampacchia [27, p. 66]. The regularity result of De Giorgi says that any $u$ that satisfies

$$\int_{\mathbb{R}^d} A(x) Du \cdot D\varphi(x) \, dx = \int_{\mathbb{R}^d} G \cdot D\varphi \, dx$$

for all $\varphi \in C_c^\infty(\Omega)$ is necessarily locally Hölder continuous. The main idea underlying the proof, and relevant to our considerations, is that one should test the equation with $\varphi = (u - k)^+$ and show that it decreases the energy, i.e.

$$\int_{\mathbb{R}^d} A(x) D\varphi \cdot D\varphi(x) \, dx \leq \int_{\mathbb{R}^d} A(x) Du \cdot D\varphi(x) \, dx.$$

This inequality allows one to leverage the classical Sobolev inequality against the equation, a 'reverse' Sobolev inequaliy, which produces the desired result, first that the solution is bounded and then that it is Hölder continuous.

Our interest in this lowering of energy with respect to such test functions stems from the relationship of the bilinear operators in question and the notion of Dirichlet forms, an idea we now explore. For any tensor $A : \mathbb{R}^d \to \mathbb{R}^{d \times d}$, define the map $B_s : H^s(\mathbb{R}^d) \times H^s(\mathbb{R}^d) \to \mathbb{R}$ by

$$B_s[u, v] := \int_{\mathbb{R}^d} A(x) D^s u \cdot D^s v(x) \, dx.$$

Then $B_s$ is a bilinear form on the Hilbert space $H^s(\mathbb{R}^d)$, while if we additionally assume $A$ is symmetric, bounded and elliptic (satisfies the lower bound in (6)), then $B_s$ satisfies

1. $B_s[u, u] \geq 0$ for all $u \in H^s(\mathbb{R}^d)$.

2. $B_s[u, v] \leq C\|D^s u\|_{L^2(\mathbb{R}^d;\mathbb{R}^d)}\|D^s v\|_{L^2(\mathbb{R}^d;\mathbb{R}^d)}$.

3. $\{u \in H^s(\mathbb{R}^d) : u \equiv 0 \text{ in } \Omega^c\}$ equipped with the scalar product

$$< u, v > := B_s[u, v]$$

   is a Hilbert space.

These are three of the four defining conditions of a Dirichlet form [21, p.4-5], and therefore one wonders

**Open Problem 1.8.** If we define $v := \min\{\max\{u, 0\}, 1\}$, can one show

$$B_s[v, v] \leq B_s[u, u]?$$

An affirmative answer to this question would imply that $B_s$ is indeed a Dirichlet form and so, in particular, by a formula of Beurling and Deny (see, for example, [21, p. 48,51]), one obtains the representation

$$
\begin{aligned}
B_s[u, v] = &\sum_{i,j=1}^{d} \int_{\mathbb{R}^d} \frac{\partial u}{\partial x_i} \frac{\partial v}{\partial x_j} d\nu_{ij}(x) \\
&+ \int_{\mathbb{R}^d} \int_{\mathbb{R}^d} (u(x) - u(y))(v(x) - v(y)) \, dJ(x, y) + \int_{\mathbb{R}^d} u(x)v(x) \, dk(x),
\end{aligned}
$$

for Radon measures $J, k, \nu_{ij}$ satisfying $J$ symmetric and positive off the diagonal, $k$ positive, and $\nu_{ij}$ such that for any compact set $K$

$$\sum_{i,j=1}^{d} \nu_{ij}(K)\xi_i\xi_j \geq 0 \quad \nu_{ij} = \nu_{ji}.$$

When $\{\nu_{ij}\}$ are absolutely continuous with respect to the Lebesgue measure and the coercivity is non-degenerate, the first term on the right hand side of the Beurling-Deny represtation falls into the framework of the De Giorgi regularity theory, while a non-local analogue to this theory has been developed for the second term by Kassmann [25, 26]. Precisely, if $dJ(x, y) = k(x, y)dxdy$ is a locally integrable kernel that satisfies

$$
\begin{aligned}
&k(x, y) = k(y, x) \\
&\lambda \le k(x, y)|x - y|^{-d-\alpha} \le \Lambda \quad |x - y| \le 1 \\
&k(x, y) \le M|x - y|^{-d-\eta} \quad |x - y| > 1
\end{aligned}
\tag{13}
$$

for some $\alpha \in (0, 2)$, $0 < \lambda \le \Lambda < +\infty$, and $\eta > 0$ and $u$ satisfies

$$
\int_{\mathbb{R}^d} \int_{\mathbb{R}^d} (u(x) - u(y))(\varphi(x) - \varphi(y))k(x, y) \, dy dx = \int_\Omega f\varphi \, dx
$$

for $f$ sufficiently nice and all $\varphi \in C_c^\infty(\Omega)$, then $u \in C_{loc}^\alpha(\Omega)$, the $L_{loc}^\infty$ estimates having been obtained previously by Fukushima [20].

A more direct approach that would immediately yield regularity is given in

**Open Problem 1.9.** Given $A : \mathbb{R}^d \to \mathbb{R}^{d \times d}$, bounded, measurable and elliptic (satisfying (6)), can one find $k_A$ that satisfies (13) such that

$$
\int_{\mathbb{R}^d} A(x)D^s u \cdot D^s \varphi(x) \, dx = \int_{\mathbb{R}^d} \int_{\mathbb{R}^d} (u(x) - u(y))(\varphi(x) - \varphi(y))k_A(x, y) \, dy dx
$$

for all $u, \varphi \in C_c^\infty(\Omega)$?

Even a negative answer to Open Problems 1.8 and 1.9 would be interesting, in that it would give a family of examples of bilinear forms whose equations exhibit regularity properties (in the case $A \in VMO$ and $p$ sufficiently large, for example) that are not Dirichlet forms.

### 1.3  Integral Functionals of the Fractional Gradient

Finally we here broaden the existence theory established in [42] to the case of possibly non-linear dependence in the field:

$$
F_s(u) = \int_{\mathbb{R}^d} f(x, u, D^s u) \, dx.
\tag{14}
$$

To this end we require a lower semicontinuity result for functionals with respect to strong-weak convergence on unbounded domains. In principle, we would like to apply Theorem 7.5 on p. 492 in [18]. However, with the introduction of an unbounded domain, one finds the constant functions are no longer integrable, and so the reduction to the case where the integrand is bounded below cannot be applied directly. If one assumes the integrand is non-negative, then the argument can be copied verbatim. We prefer to keep the general assumptions of the theorem, supplementing them with the simple additional assumption that outside a large ball, the integrand $f$ has the lower bound

$$f(x, z, \xi) \geq \alpha(x) + \tilde{\beta}(x) \cdot \xi - C|z|^q,$$

for some $\alpha \in L^1(\mathbb{R}^d)$, $\tilde{\beta} \in L^{p/(p-1)}(\mathbb{R}^d)$ and $C > 0$. This only differs from the standard theorem in that $\tilde{\beta}$ does not depend on $z$. In particular, we prove the following

**Theorem 1.10.** *Let* $1 < p, q < \infty$ *and suppose that* $f : \mathbb{R}^d \times \mathbb{R} \times \mathbb{R}^d \to \mathbb{R}$ *is* $\mathcal{L}^d \times \mathcal{B}$ *measurable function such that*

$$f(x, z, \xi) \geq -C(|z|^q + |\xi|^p) - \omega(x)$$

*for* $\mathcal{L}^d$ *almost every* $x$ *and for all* $(z, \xi) \in \mathbb{R} \times \mathbb{R}^d$, *for some* $\omega \in L^1(\mathbb{R}^d)$ *and* $C > 0$. *Assume that* $f(x, \cdot, \cdot)$ *is lower semicontinuous in* $\mathbb{R} \times \mathbb{R}^d$ *for* $\mathcal{L}^d$ *almost every* $x \in \mathbb{R}^d$. *Further assume that*

    *i. $f(x, z, \cdot)$ is convex in $\mathbb{R}^d$ for $\mathcal{L}^d$ almost every $x \in \mathbb{R}^d$ and for all $z \in \mathbb{R}$;*

    *ii. For $\mathcal{L}^d$ almost every $x \in \mathbb{R}^d$ and all $(z, \xi) \in \mathbb{R} \times \mathbb{R}^d$,*

$$f(x, z, \xi) \geq \alpha(x) - C|z|^q + \beta(x, z) \cdot \xi$$

    *where $\alpha \in L^1(\mathbb{R}^d)$, $\beta : \mathbb{R}^d \times \mathbb{R} \to \mathbb{R}^d$ is a $\mathcal{L}^d \times \mathcal{B}$ measurable function such that $\beta(x, z) \equiv \tilde{\beta}(x)$ outside $B(0, R)$ for some $R > 0$ with $\tilde{\beta} \in L^{p/(p-1)}(B(0, R)^c)$ and $C > 0$;*

    *iii. For any sequences $\{u_n\} \subset L^q(\mathbb{R}^d)$ $\{v_n\} \subset L^p(\mathbb{R}^d; \mathbb{R}^d)$ such that $u_n \to u$ strongly in $L^q(\mathbb{R}^d)$ and $v_n \to v$ weakly in $L^p(\mathbb{R}^d; \mathbb{R}^d)$, and such that*

$$\sup_n \int_{B(0,R)} f(x, u_n(x), v_n(x)) \, dx < +\infty,$$

    *then the sequence $|\beta(x, u_n(x))|^{p/(p-1)}$ is equi-integrable in $B(0, R)$.*

*Then the functional*

$$(u, v) \in L^q(\mathbb{R}^d) \times L^p(\mathbb{R}^d; \mathbb{R}^d) \mapsto \int_{\mathbb{R}^d} f(x, u, v) \, dx$$

*is sequentially lower semicontinuous with respect to strong convergence in $L^q(\mathbb{R}^d)$ and weak convergence in $L^p(\mathbb{R}^d; \mathbb{R}^d)$.*

Let us next recall the definition of the spaces on which we will show the existence of minimizers. Here and in what follows, $\Omega \subset \mathbb{R}^d$ is a bounded open set. We first introduce the fractional Sobolev spaces without boundary values

$$H^{s,p}(\mathbb{R}^d) := \{u \in L^p(\mathbb{R}^d) : D^s u \in L^p(\mathbb{R}^d; \mathbb{R}^d)\},$$

and those with a given boundary value $g \in H^{s,p}(\mathbb{R}^d)$

$$H_g^{s,p}(\Omega) := \{u \in H^{s,p}(\mathbb{R}^d) : u \equiv g \text{ in } \Omega^c\}.$$

**Remark 1.11.** We have changed notation here in contrast to the preceding paper, where $L^{s,p}(\mathbb{R}^d)$ and $L_g^{s,p}(\Omega)$ were used to denote the two previous spaces, respectively.

Then the next result of this section is the following theorem on the existence of minimizers to the general integral functionals (14):

**Theorem 1.12.** *Assume $f : \mathbb{R}^d \times \mathbb{R} \times \mathbb{R}^d$ is $\mathcal{L}^d \times \mathcal{B}$ measurable, $f(x, \cdot, \cdot)$ is lower semicontinuous for $\mathcal{L}^d$ almost every $x \in \mathbb{R}^d$, and $f(x, z, \cdot)$ is convex for $\mathcal{L}^d$ almost every $x \in \mathbb{R}^d$ and all $z \in \mathbb{R}$. Further assume that $f$ satisfies the coercivity condition*

$$f(x, z, \xi) \geq a(x) + b|\xi|^p$$

*for almost every $x \in \mathbb{R}^d$, for every $(z, \xi) \in \mathbb{R} \times \mathbb{R}^d$, and for some $a \in L^1(\mathbb{R}^d)$, $b > 0$, and $p > 1$. Finally, assume that $F_s(u_0)$ is finite for some $u_0 \in H_g^{s,p}(\Omega)$. Then there exist at least one minimizer $u \in H_g^{s,p}(\Omega)$ of the functional $F_s$:*

$$F_s(u) \leq F_s(v)$$

*for all $v \in H_g^{s,p}(\Omega)$.*

Accordingly, we obtain existence of solutions to corresponding Euler-Lagrange equations under further smoothness and growth assumptions on $f$.

**Theorem 1.13.** *Assume that the function $f : \mathbb{R}^d \times \mathbb{R} \times \mathbb{R}^d$ satisfies the conditions of Theorem 1.12. Additionally assume growth conditions*

$$\begin{aligned}
|f(x, z, \xi)| &\leq C(|z|^p + |\xi|^p) + \gamma_1(x), \\
|D_z f(x, z, \xi)| &\leq C(|z|^{p-1} + |\xi|^{p-1}) + \gamma_2(x), \\
|D_\xi f(x, z, \xi)| &\leq C(|z|^{p-1} + |\xi|^{p-1}) + \gamma_3(x).
\end{aligned} \tag{15}$$

where $\gamma_1, \gamma_2, \gamma_3 \in L^{p/(p-1)}(\mathbb{R}^d)$. *Further suppose that $u$ is a minimizer of $F_s$ over $\in H_g^{s,p}(\Omega)$. Then $u$ satisfies*

$$\int_{\mathbb{R}^d} f_z(x, u, D^s u)\, \varphi + D_\xi f(x, u, D^s u) \cdot D^s \varphi \, dx = 0$$

*for all $\varphi \in C_c^\infty(\Omega)$.*

The plan of the papers is as follows. We first define some notation and prove an important tool - a compactness result - in Section 2. We then prove the Hardy inequality in Section 3 followed by the regularity results in Section 4. Finally we give the proofs of our results concerning integral functionals of the fractional gradient in Section 5.

## 2 Preliminaries

In this paper we work in $d$-dimensional Euclidean space, denoting by $\mathcal{L}^d$ the Lebesgue measure, which we often shorten to $dx$ in the integration formulas.
We denote by $\mathcal{B}$ the Borel $\sigma$-algebra.
We write $B(x, r)$ for a ball centered at $x$ with radius $r > 0$.
In the introduction we have utilized the constants $c_{d,s}, \tilde{c}_{d,s}$ to ensure that

$$(D^s u)^\wedge(\xi) = -2\pi i \xi (2\pi|\xi|)^{-1+s} \hat{u}(\xi),$$

$$((-\Delta)^{s/2} u)^\wedge(\xi) = (2\pi|\xi|)^s \hat{u}(\xi),$$

where we use the convention

$$\widehat{u}(\xi) = \int_{\mathbb{R}^N} u(x) e^{-2\pi i x \cdot \xi} \, dx$$

In particular, one has

$$\tilde{c}_{d,s} := \frac{2^{s-1} s \Gamma(\frac{d+s}{2})}{\pi^{\frac{d}{2}} \Gamma(1 - s/2)},$$

$$c_{d,s} := (-d + 1 - s) \frac{1}{\gamma(1-s)}.$$

Here, the constant $\gamma$ arises from the consideration of $\alpha < 0$, where the fractional Laplacian has as its inverse the Riesz potential $I_{-\alpha} u := (-\Delta)^{\alpha/2} u$, which has integral formula for $s \in (0, d)$

$$I_\alpha u(x) = \frac{1}{\gamma(\alpha)} \int_{\mathbb{R}^d} \frac{u(y)}{|x-y|^{d-\alpha}} \, dy,$$

and

$$\gamma(\alpha) = \frac{\pi^{d/2} 2^\alpha \Gamma(\alpha/2)}{\Gamma(\frac{d-\alpha}{2})}.$$

An important step in the argument that one has existence of minimizers of integral functionals with non-linear dependence in the field with respect to weak convergence in $H_g^{s,p}(\Omega)$ is the following compactness result that improves weak convergence in the fields to strong convergence.

**Theorem 2.1.** *Assume $\Omega$ is a bounded open subset of $\mathbb{R}^d$, suppose $s \in (0,1)$ and $1 < p < \frac{d}{s}$. Then for any sequence $\{u_n\} \subset H_g^{s,p}(\Omega)$ such that*

$$u_n \to u \text{ weakly in } H_g^{s,p}(\Omega),$$

*we have that*

$$u_n \to u \text{ strongly in } L^q(\Omega)$$

*for every $q \in [1, p^*)$. Here $\frac{1}{p^*} = \frac{1}{p} - \frac{s}{d}$.*

*Proof of Theorem* 2.1. Suppose $u_m \to u$ weakly in $H_g^{s,p}(\Omega)$. We will show that for any subsequence (which we will not relabel), there is a further subsequence which is Cauchy in $L^q(\Omega)$, and therefore the original sequence is strongly convergent. Without loss of generality, we replace the sequence $u_m - g$ with $u_m \in H_0^{s,p}(\Omega)$. Let $\eta_\epsilon$ be a standard mollifier. Set

$$u_m^\epsilon = \eta_\epsilon * u_m$$

for $\epsilon > 0$ and $m \in \mathbb{N}$. We may assume all functions $\{u_m\}_{m=1}^\infty$ have support in $\Omega'$. We claim that

$$u_m^\epsilon \to u_m \quad \text{uniformly in } L^q(\Omega) \text{ as } \epsilon \to 0.$$

By density of smooth functions with rapidly decreasing decay at infinity, we can represent $u_m = I_s(-\Delta)^{s/2} u_m$. Define $v_m = (-\Delta)^{s/2} u_m$. Then boundedness of the Riesz transforms on $L^p(\mathbb{R}^d)$ for $1 < p < +\infty$, we have that $v_m$ and $|D^s u_m|$ have comparable norms in $L^p(\mathbb{R}^d)$. Thus, we estimate

$$u_m^\epsilon(x) - u_m(x) = \int_{B(0,1)} \eta(y)(u_m(x - \epsilon y) - u_m(x)) \, dy$$

$$= \int_{B(0,1)} \eta(y) \int_{\mathbb{R}^d} \frac{v_m(x - \epsilon y - z)}{|z|^{d-s}} - \frac{v_m(x-z)}{|z|^{d-s}} \, dz \, dy$$

$$= \int_{B(0,1)} \eta(y) \int_{\mathbb{R}^d} \left( \frac{1}{|z - \epsilon y|^{d-s}} - \frac{1}{|z|^{d-s}} \right) v_m(x-z) \, dz \, dy.$$

Changing variables $z = \epsilon \tilde{z}$, we find

$$
\begin{aligned}
u_m^\epsilon(x) - u_m(x) &= \epsilon^d \int_{B(0,1)} \eta(y) \int_{\mathbb{R}^d} \left( \frac{1}{|\epsilon \tilde{z} - \epsilon y|^{d-s}} - \frac{1}{|\epsilon \tilde{z}|^{d-s}} \right) v_m(x - z) \, d\tilde{z} \, dy \\
&= \frac{\epsilon^d}{\epsilon^{d-s}} \int_{B(0,1)} \eta(y) \int_{\mathbb{R}^d} \left( \frac{1}{|\tilde{z} - y|^{d-s}} - \frac{1}{|\tilde{z}|^{d-s}} \right) v_m(x - \epsilon \tilde{z}) \, d\tilde{z} \, dy \\
&= \epsilon^s \int_{B(0,1)} \eta(y) \int_{\mathbb{R}^d} \left( \frac{1}{|\tilde{z} - y|^{d-s}} - \frac{1}{|\tilde{z}|^{d-s}} \right) v_m(x - \epsilon \tilde{z}) \, d\tilde{z} \, dy.
\end{aligned}
$$

Thus, integrating over a bounded open set $\Omega'$ which contains $\Omega$ and we obtain

$$
\begin{aligned}
&\int_{\Omega'} |u_m^\epsilon(x) - u_m(x)| \, dx \\
&\leq \epsilon^s \int_{\Omega'} \int_{B(0,1)} \int_{\mathbb{R}^d} \eta(y) \left| \frac{1}{|z - y|^{d-s}} - \frac{1}{|z|^{d-s}} \right| |v_m(x - \epsilon z)| \, dz \, dy \, dx \\
&= \epsilon^s \int_{B(0,1)} \int_{\mathbb{R}^d} \eta(y) \left| \frac{1}{|z - y|^{d-s}} - \frac{1}{|z|^{d-s}} \right| \int_{\Omega'} |v_m(x - \epsilon z)| \, dx \, dz \, dy \\
&\leq \epsilon^s |\Omega'|^{1/p'} \int_{B(0,1)} \int_{\mathbb{R}^d} \eta(y) \left| \frac{1}{|z - y|^{d-s}} - \frac{1}{|z|^{d-s}} \right| \left( \int_{\Omega'} |v_m(x - \epsilon z)|^p \, dx \right)^{1/p} \, dz \, dy \\
&\leq \epsilon^s |\Omega'|^{1/p'} \int_{B(0,1)} \int_{\mathbb{R}^d} \eta(y) \left| \frac{1}{|z - y|^{d-s}} - \frac{1}{|z|^{d-s}} \right| \, dz \, dy \left( \int_{\mathbb{R}^d} |v_m(x)|^p \, dx \right)^{1/p} \\
&\leq C \epsilon^s |\Omega'|^{1/p'} \left( \sup_{y \in B(0,1)} \int_{\mathbb{R}^d} \left| \frac{1}{|z - y|^{d-s}} - \frac{1}{|z|^{d-s}} \right| \, dz \right) \|v_m\|_{L^p(\mathbb{R}^d)}
\end{aligned}
$$

For $y \in B(0,1)$, we estimate the following two integrals

$$
\int_{B(0,2)} \left| \frac{1}{|z - y|^{d-s}} - \frac{1}{|z|^{d-s}} \right| \, dz \leq 2 \int_{B(0,2)} \frac{1}{|z|^{d-s}} \, dz < +\infty
$$

and

$$
\begin{aligned}
\int_{\mathbb{R}^d \setminus B(0,2)} \left| \frac{1}{|z-y|^{d-s}} - \frac{1}{|z|^{d-s}} \right| dz &= \int_{\mathbb{R}^d \setminus B(0,2)} \left| \int_0^1 \frac{d}{dt} \frac{1}{|z-ty|^{d-s}} \, dt \right| dz \\
&= \int_{\mathbb{R}^d \setminus B(0,2)} \left| \int_0^1 \frac{(ty-z) \cdot y}{|z-ty|^{s+2-s}} \, dt \right| dz \\
&\leq \int_{\mathbb{R}^d \setminus B(0,2)} \int_0^1 \frac{|y|}{|z-ty|^{d+1-s}} \, dt \, dz \\
&\leq |y| \int_0^1 \int_{\mathbb{R}^d \setminus B(0,2)} \frac{1}{|z-ty|^{d+1-s}} \, dz \, dt \\
&\leq |y| \int_0^1 \int_{\mathbb{R}^d \setminus B(0,1)} \frac{1}{|z|^{d+1-s}} \, dz \, dt < +\infty.
\end{aligned}
$$

Since $\{u_m\}$ is weakly convergent, we know that it is a bounded sequence in $H_g^{s,p}(\Omega)$ and so $\|v_m\|_{L^p(\mathbb{R}^d)}$ is bounded. Thus, we find that

$$
\|u_m^\epsilon - u_m\|_{L^1(\Omega)} = O(\epsilon^s).
$$

On the other hand, the Sobolev inequality says that

$$
\begin{aligned}
\|u_m^\epsilon - u_m\|_{L^{p^*}(\Omega)} &\leq \|u_m^\epsilon - u_m\|_{L^{p^*}(\mathbb{R}^d)} \\
&\leq C\|D^s(u_m^\epsilon - u_m)\|_{L^p(\mathbb{R}^d)} \leq C\|D^s u_m\|_{L^p(\mathbb{R}^d)} < +\infty.
\end{aligned}
$$

The previous $L^1(\Omega)$ bound and the interpolation inequality

$$
\|u_m^\epsilon - u_m\|_{L^q(\Omega)} \leq \|u_m^\epsilon - u_m\|_{L^1(\Omega)}^\theta \|u_m^\epsilon - u_m\|_{L^{p^*}(\Omega)}^{1-\theta}
$$

implies that one has, for any $1 \leq q < p^*$

$$
\|u_m^\epsilon - u_m\|_{L^q(\Omega)} \leq C\epsilon^{s\theta}
$$

where the constant $C$ is independent of $m$. Here, precisely $\frac{1}{q} = \theta + \frac{1-\theta}{p^*}$ for some $0 < \theta < 1$.

We would now like to invoke the Arzela-Ascoli theorem concerning the sequence $\{u_m^\epsilon\}$ of smooth functions restricted on $\overline{\Omega}$ for every fixed $\epsilon$. We therefore prove that for each fixed $\epsilon > 0$, the sequence $\{u_m^\epsilon\}$ is uniformly bounded and equicontinuous on $\overline{\Omega}$. For $x \in \Omega$, we estimate

$$|u_m^\epsilon(x)| \leq \int_{B(x,\epsilon)} \eta_\epsilon(x-y)|u_m(y)|\,dy$$

$$\leq \|\eta_\epsilon\|_{L^\infty(\mathbb{R}^d)}\|u_m\|_{L^1(\Omega')} \leq \frac{C}{\epsilon^d} < +\infty$$

for $m \in \mathbb{N}$. Moreover, since $\eta_\epsilon$ are smooth and $u_m$ have compact support, we have

$$\nabla u_m^\epsilon = \int_{\mathbb{R}^d} \nabla \eta_\epsilon(x-y)u_m(y)\,dy$$

and therefore

$$|\nabla u_m^\epsilon| \leq \|\nabla \eta_\epsilon\|_{L^\infty(\mathbb{R}^d)}\|u_m\|_{L^1(\Omega)}$$

for $m \in \mathbb{N}$. These estimates prove the claim of the uniformly boundedness and equicontinuity of the sequence $\{u_m^\epsilon\}_{m=1}^\infty$ for every fixed $\epsilon$.

In the final step, we want to construct a subsequence $\{u_{m_k}\}_{k=1}^\infty \subset \{u_m\}_{m=1}^\infty$ such that

$$\limsup_{j,k\to\infty} \|u_{m_j} - u_{m_k}\|_{L^q(\Omega')} = 0.$$

In order to show this, first we claim for fixed $\delta$, there exists subsequence $\{u_{m_k}\}_{k=1}^\infty \subset \{u_m\}_{m=1}^\infty$ such that

$$\limsup_{j,k\to\infty} \|u_{m_j} - u_{m_k}\|_{L^q(\Omega')} \leq \delta.$$

For $\epsilon$ small enough, we have

$$\|u_{m_k}^\epsilon - u_{m_k}\|_{L^q(\Omega')} \leq \frac{\delta}{2}.$$

Since $\{u_m^\epsilon\}$ have support in some fixed bounded set $\Omega'$, we apply the Arzela-Ascoli theorem to find a subsequence $\{u_{m_k}^\epsilon\}_{k=1}^\infty \subset \{u_m^\epsilon\}_{m=1}^\infty$ converging uniformly in $\Omega'$. This is

$$\limsup_{j,k\to\infty} \|u_{m_j}^\epsilon - u_{m_k}^\epsilon\|_{L^q(\Omega')} = 0.$$

Therefore, we have

$$
\begin{aligned}
\limsup_{j,k\to\infty} & \|u_{m_j} - u_{m_k}\|_{L^q(\Omega')} \\
& \leq \limsup_{j\to\infty} \|u_{m_j}^\epsilon - u_{m_j}\|_{L^q(\Omega')} + \limsup_{j,k\to\infty} \|u_{m_j}^\epsilon - u_{m_k}^\epsilon\|_{L^q(\Omega')} \\
& \quad + \limsup_{k\to\infty} \|u_{m_k}^\epsilon - u_{m_k}\|_{L^q(\Omega')} \\
& \leq \frac{\delta}{2} + \frac{\delta}{2} = \delta.
\end{aligned}
$$

Thus, choosing the sequence $\delta_n := \frac{1}{n}$ and a standard diagonalization argument, we may find a subsequence $\{u_{m_k}\}_{k=1}^\infty \subset \{u_m\}_{m=1}^\infty$ such that

$$
\limsup_{j,k\to\infty} \|u_{m_j} - u_{m_k}\|_{L^q(\Omega')} = 0.
$$

This shows the sequence is Cauchy, which by completeness of $L^q(\Omega)$ implies the strong convergence of the sequence to some function, which by uniqueness of the weak limit implies $u_{m_j} \to u$ strongly in $L^q(\Omega)$ (and also all of $\mathbb{R}^d$, since $u_m = u \equiv 0$ in $\Omega^c$). $\qquad\square$

## 3   Hardy's Inequality

*Proof of Theorem* 1.2.  Letting $C_{d,s}$ to denote the constant on the left hand side, we first show that

$$
C_{d,s} \int_{\mathbb{R}^d} \frac{|u|}{|x|^s} \, dx = \int_{\mathbb{R}^d} -D^s|u| \cdot \frac{x}{|x|} \, dx,
$$

from which the inequality follows by bringing the modulus into the integral. We have

$$
\begin{aligned}
C_{d,s} \int_{\mathbb{R}^d} \frac{|u|}{|x|^s} \, dx &= C_{d,s}\gamma(d-s) \int_{\mathbb{R}^d} |u| \cdot I_{d-s} \, dx \\
&= C_{d,s}\gamma(d-s) \int_{\mathbb{R}^d} I_{1-s} * |u| \cdot I_{d-1} \, dx \\
&= C_{d,s}\frac{\gamma(d-s)}{\gamma(d-1)} \int_{\mathbb{R}^d} I_{1-s} * |u| \cdot \frac{1}{|x|} \, dx.
\end{aligned}
$$

Now recalling that

$$
div\,\frac{x}{|x|} = (d-1)\frac{1}{|x|},
$$

we find

$$C_{d,s} \frac{\gamma(d-s)}{\gamma(d-1)} \int_{\mathbb{R}^d} I_{1-s} * |u| \cdot \frac{1}{|x|} \, dx$$

$$= \frac{1}{(d-1)} C_{d,s} \frac{\gamma(d-s)}{\gamma(d-1)} \int_{\mathbb{R}^d} I_{1-s} * |u| \cdot div \frac{x}{|x|} \, dx$$

$$= \frac{1}{(d-1)} C_{d,s} \frac{\gamma(d-s)}{\gamma(d-1)} \int_{\mathbb{R}^d} -D I_{1-s} * |u| \cdot \frac{x}{|x|} \, dx$$

$$= \frac{1}{(d-1)} C_{d,s} \frac{\gamma(d-s)}{\gamma(d-1)} \int_{\mathbb{R}^d} -D^s |u| \cdot \frac{x}{|x|} \, dx,$$

and the claim is proven since the constant $C_{d,s}$ is defined such that the coefficient of the right hand side is one.                                                  □

## 4   Regularity

The following fundamental result underlies the regularity of homogeneous fractional equations for which the regularity is known in the corresponding non-fractional setting.

**Proposition 4.1.** *Let* $\Omega_1 \subset\subset \Omega_2 \subset\subset \Omega$, $\phi \in C_c^\infty(\Omega_1), \eta \in C_c^\infty(\Omega)$ *be such that* $\eta \equiv 1$ *on* $\Omega_2$. *Then denoting by* $\eta^c := (1 - \eta)$, *the operator* $T$ *defined by*

$$T(\phi) := D^s(\eta^c(-\Delta)^{\frac{1-s}{2}} \phi)$$

*is bounded from functions* $\phi \in L^1(\mathbb{R}^d)$ *with* $\operatorname{supp} \phi \subset \Omega_1$ *into* $L^p(\mathbb{R}^d; \mathbb{R}^d)$. *In particular, one has the bound*

$$\|T(\phi)\|_{L^p(\mathbb{R}^d;\mathbb{R}^d)} \le C_{\Omega_1,\Omega_2,d,s,p,\eta} \|\phi\|_{L^1(\Omega_1)}. \tag{16}$$

This proposition has been established in the paper [41], whose argument we repeat here for the convenience of the reader.

*Proof.* We will show that for $T$ as defined above, one has the estimate

$$\|T(\phi)\|_{L^q(\mathbb{R}^d;\mathbb{R}^d)} \le C_{\Omega_1,\Omega_2,d,s,p} \|\phi\|_{L^1(\mathbb{R}^d)}. \tag{17}$$

We use the disjoint support arguments as in [4, Lemma A.1] [32, Lemma 3.6.]: First we see that since $\eta^c(x)\phi(x) \equiv 0$,

$$T(\phi) = \tilde{c}_{d,1-s} D^s \int_{\mathbb{R}^d} \frac{-\eta^c(x)\phi(y)}{|x-y|^{d+1-s}} \, dy.$$

Now taking a cutoff-function $\zeta$ whose support has a positive distance from the boundary of $\Omega_2$, $\zeta \equiv 1$ on $\Omega_1$ we have

$$T(\phi) = \tilde{c}_{d,1-s} \, D^s \int_{\mathbb{R}^d} \frac{-\eta^c(x)\zeta(y)\phi(y)}{|x-y|^{d+1-s}} \, dy = \tilde{c}_{d,1-s} \int_{\mathbb{R}^d} k(x,y) \, \phi(y) \, dy,$$

where

$$\kappa(x,y) := \frac{-\eta^c(x)\,\zeta(y)}{|x-y|^{d+1-s}} \quad \text{and} \quad k(x,y) := D_x^s \kappa(x,y).$$

The positive distance between the supports of $\eta^c$ and $\zeta$ implies that these kernels $k$, $\kappa$ are a smooth, bounded, integrable (both, in $x$ and in $y$), and thus by a Young-type convolution argument we obtain (17). One can also argue by interpolation, as Minkowski's inequality for integrals implies

$$\left\| \int_{\mathbb{R}^d} D_x^s \kappa(\cdot, y) \, \phi(y) \, dy \right\|_{L^p(\mathbb{R}^d; \mathbb{R}^d)} \leq \sup_y \| D_x^s \kappa(\cdot, y) \|_{L^p(\mathbb{R}^d; \mathbb{R}^d)} \| \phi \|_{L^1(\mathbb{R}^d)},$$

while Theorem 2.4 on p. 886 of Adams and Meyers paper [1] can be applied to obtain

$$\| D_x^s \kappa(\cdot, y) \|_{L^p(\mathbb{R}^d; \mathbb{R}^d)} \leq C \, \| D_x \kappa(\cdot, y) \|_{L^p(\mathbb{R}^d; \mathbb{R}^d)}^s \, \| R\kappa(\cdot, y) \|_{L^p(\mathbb{R}^d; \mathbb{R}^d)}^{1-s}.$$

Then boundedness of the Riesz Transform and integrability of the kernels establishes (17) and the proof is finished.                                                   □

As a consequence, we deduce

**Corollary 4.2.** *Let $\Omega_1 \subset\subset \Omega_2 \subset\subset \Omega$, $\phi \in C_c^\infty(\Omega_1), \eta \in C_c^\infty(\Omega)$ be such that $\eta \equiv 1$ on $\Omega_2$. Then denoting by $T^*$ the adjoint operator to*

$$T(\phi) := D^s(\eta^c(-\Delta)^{\frac{1-s}{2}}\phi),$$

*one has*

$$T^* : L^q(\mathbb{R}^d; \mathbb{R}^d) \to L^\infty(\Omega_1)$$

*for every $1 < q < +\infty$, with the operator norm of $T^*$ depending on $\Omega_1, \Omega_2, d, s, q$ and $\eta$.*

*Proof of Theorem 1.5.* Suppose that $u \in H_g^{s,p}(\Omega)$ satisfies

$$\int_{\mathbb{R}^d} a(x, D^s u) \cdot D^s \varphi = 0 \tag{18}$$

for all $\varphi \in C_c^\infty(\Omega)$. Define $v := I_{1-s}u$, where $I_\sigma$ is the Riesz potential, the inverse of $(-\Delta)^{\sigma/2}$. Now let $\Omega_1 \subset \Omega$ be an arbitrary open set compactly contained in $\Omega$, and let $\phi$ be a test function supported in $\Omega_1$. Pick an open set $\Omega_2$ so that $\Omega_1 \subset \Omega_2 \subset \Omega$ and a cutoff function $\eta$, supported in $\Omega$ and constantly one in $\Omega_2$. Then in particular one can take

$$\varphi := \eta(-\Delta)^{\frac{1-s}{2}}\phi$$

as a test function to obtain

$$\int_{\mathbb{R}^d} a(x, D^s v) \cdot D^s(\eta(-\Delta)^{\frac{1-s}{2}}\phi) = 0.$$

Thus,

$$\int_{\mathbb{R}^d} a(x, D^s v) \cdot D\phi = \int_{\mathbb{R}^d} a(x, D^s v) \cdot D^s(\eta^c(-\Delta)^{\frac{1-s}{2}}\phi)$$

where $\eta^c := (1 - \eta)$. We set

$$T(\phi) := D^s(\eta^c(-\Delta)^{\frac{1-s}{2}}\phi),$$

and from the assumptions on $a$ we may apply Corollary 4.2 to deduce that

$$T^* \cdot a(x, D^s u) \in L_{loc}^\infty(\Omega).$$

In other words, $v$ is a solution to the equation

$$\int_\Omega a(x, Dv) \cdot D\phi \, dx = \int_\Omega \phi \, d\mu.$$

Thus by Theorems 1.4 and 1.6 in [29], we find that $\mu \in L_{loc}^\infty(\Omega)$ implies $Dv \in C_{loc}^\alpha(\Omega)$. Now as $Dv = D^s u \in C_{loc}^\alpha(\Omega)$, we obtain $u \in C_{loc}^{s+\alpha}(\Omega)$. □

*Proof of Theorem* 1.6. Suppose that $u \in H_g^{s,2}(\Omega)$ satisfies

$$\int_{\mathbb{R}^d} A(x)D^s u \cdot D^s\varphi = \int_{\mathbb{R}^d} G \cdot D^s\varphi \tag{19}$$

for all $\varphi \in C_c^\infty(\Omega)$. Define $v := I_{1-s}u$, where $I_\sigma$ is the Riesz potential, the inverse of $(-\Delta)^{\sigma/2}$. Now let $\Omega_1 \subset \Omega$ be an arbitrary open set compactly contained in $\Omega$, and let $\phi$ be a test function supported in $\Omega_1$. Pick an open set $\Omega_2$ so that $\Omega_1 \subset \Omega_2 \subset \Omega$ and a cutoff function $\eta$, supported in $\Omega$ and constantly one in $\Omega_2$. Then in particular one can take

$$\varphi := \eta(-\Delta)^{\frac{1-s}{2}}\phi$$

as a test function in (19) to obtain

$$\int_{\mathbb{R}^d} A(x)Dv \cdot D^s(\eta(-\Delta)^{\frac{1-s}{2}}\phi) = \int_{\mathbb{R}^d} G \cdot D^s(\eta(-\Delta)^{\frac{1-s}{2}}\phi)$$

That is,

$$\int_{\mathbb{R}^d} A(x)Dv \cdot D\phi = \int_{\mathbb{R}^d} A(x)Dv \cdot D^s(\eta^c(-\Delta)^{\frac{1-s}{2}}\phi) + G \cdot D^s(\eta(-\Delta)^{\frac{1-s}{2}}\phi)$$

$$= \int_{\mathbb{R}^d} G \cdot D\phi + (A(x)Dv - G) \cdot D^s(\eta^c(-\Delta)^{\frac{1-s}{2}}\phi)$$

where $\eta^c := (1 - \eta)$. We set

$$T(\phi) := D^s(\eta^c(-\Delta)^{\frac{1-s}{2}}\phi),$$

and by Corollary 4.2, we find that $v$ is a solution to the classical elliptic equation with bounded and measurable coefficients

$$\int_{\Omega} A(x)Dv \cdot D\phi \; dx = \int_{\Omega} G \cdot D\phi + T^* \cdot ((A(x)Dv - G))\phi.$$

Thus, by the regularity theory known for such an equation (e.g. [27, p. 66]), we find $v \in C_{loc}^{\alpha}(\Omega)$, which is to say $I_{1-s}u \in C_{loc}^{\alpha}(\Omega)$.

<div align="right">□</div>

## 5   Integral Functionals of the Fractional Gradient

In this section, we consider the variational problem

$$\inf_{u \in H_g^{s,p}(\Omega)} \int_{\mathbb{R}^d} f(x, u, D^s u) \; dx.$$

Under suitable hypothesis, we establish the existence of minimizers, while with further assumptions we show that these minimizers satisfy corresponding Euler-Lagrange equations.

We begin by proving the lower semicontinuity result for strong-weak convergence stated in Theorem 1.10.

*Proof.* Suppose $u_n \to u$ strongly in $L^q(\mathbb{R}^d)$, $v_n \to v$ weakly in $L^p(\mathbb{R}^d; \mathbb{R}^d)$. Now we may assume that

$$\liminf_{n \to \infty} \int_{\mathbb{R}^d} f(x, u_n, v_n) \; dx < +\infty$$

or else there is nothing to prove. From the assumptions of the theorem we find $R > 0$ such that

$$f(x, z, \xi) \geq \alpha(x) - C|z|^q + \tilde{\beta}(x) \cdot \xi$$

for $\mathcal{L}^d$ almost every $x \in B(0, R)^c$. Thus, we split the integrand and use super-additivity of the limit inferior to obtain

$$\liminf_{n \to \infty} \int_{\mathbb{R}^d} f(x, u_n, v_n) \, dx \geq \liminf_{n \to \infty} \int_{B(0,R)} f(x, u_n, v_n) \, dx$$

$$+ \liminf_{n \to \infty} \int_{B(0,R)^c} f(x, u_n, v_n) \, dx$$

The first term in the integrand now satisfies the hypothesis of Theorem 7.5 in [18] with the bounded domain $E = B(0, R)$ and so we find

$$\liminf_{n \to \infty} \int_{B(0,R)} f(x, u_n, v_n) \, dx \geq \int_{B(0,R)} f(x, u, v) \, dx.$$

Meanwhile, for the second term we define the perturbation of $f$

$$\tilde{f}(x, z, \xi) := f(x, z, \xi) - \alpha(x) + C|z|^q - \tilde{\beta}(x) \cdot \xi.$$

Then $\tilde{f}$ is non-negative in $B(0, R)^c$ and so the blow-up argument in Step 1 of Theorem 7.5 can be applied in the unbounded domain $B(0, R)^c$. In particular, we are in the case $p > 1$ in the appeal to Theorem 7.2 for a representation of a coercive perturbation of $\tilde{f}$ as the supremum of affine functions, which is allowed even for unbounded domains. The rest of the argument remains unchanged, since the argument is localized by the blow-up. Thus we find that

$$\liminf_{n \to \infty} \int_{B(0,R)^c} f(x, u_n, v_n) - \alpha(x) + C|u_n|^q - \tilde{\beta}(x) \cdot v_n$$

$$\geq \int_{B(0,R)^c} f(x, u, v) - \alpha(x) + C|u|^q - \tilde{\beta}(x) \cdot v,$$

which from the strong convergence of $u_n$ and the weak convergence of $v_n$ implies

$$\liminf_{n \to \infty} \int_{B(0,R)^c} f(x, u_n, v_n) \, dx \geq \int_{B(0,R)^c} f(x, u, v) \, dx.$$

Combining this with the inequality in $B(0, R)$ concludes the proof.      □

*Proof of Theorem* 1.12. Since we have assumed that there exists a function $u_0 \in H_g^{s,p}(\Omega)$ such that $F_s(u_0) < +\infty$, we may find a minimizing sequence $\{u_k\}$ such that

$$\lim_{k \to \infty} F_s(u_k) = \inf_{u \in H_g^{s,p}(\Omega)} F_s(u) =: C_s < +\infty.$$

Then the coercivity assumption implies that the fractional gradients remain on a bounded set of $L^p(\mathbb{R}^d; \mathbb{R}^d)$: For $k$ sufficiently large, we have

$$\int_{\mathbb{R}^d} |D^s u_k|^p \leq C_s + 1.$$

Now if $sp < d$, Hölder's inequality and the fractional Sobolev inequality imply that for any $1 \leq q < p^*$

$$\|u_k - g\|_{L^q(\Omega)} \leq C\|u_k - g\|_{L^{p^*}(\Omega)} \leq C\|D^s(u_k - g)\|_{L^p(\mathbb{R}^d)},$$

while if $sp = d$, the bound also holds because of local exponential integrability of $u_k - g$. Finally, if $sp > d$, the sequence $u_k - g \in L^\infty(\mathbb{R}^d)$ by Morrey's inequality (see, for example, [42]). Thus, $\{u_k\}_{k=1}^\infty$ is bounded in $H_g^{s,p}(\Omega)$. According to the weak compactness theorem, there exist subsequence $\{u_{k_j}\}_{j=1}^\infty$ and $u \in H_g^{s,p}(\Omega)$ such that $u_{k_j}$ converges strongly to $u$ in $L^p(\mathbb{R}^d)$ and $D^s u_{k_j}$ converges weakly to $D^s u$ in $L^p(\mathbb{R}^d; \mathbb{R}^d)$. By subtracting the function $a$ in the lower bound for $f$, we find that $f$ is non-negative and satisfies the hypothesis of Theorem 1.10. In particular, taking $q = p$ we have that the functional $F_s$ is lower semicontinuous with respect to this strong-weak convergence, and so we obtain

$$F_s(u) \leq \liminf_{j \to \infty} F(u_{k_j}) = \lim_{j \to \infty} F(u_{k_j}) = \inf_{v \in H_g^{s,p}(\Omega)} F(v).$$

This shows that $u \in H_g^{s,p}(\Omega)$ minimizes the functional $F_s$.                    □

Finally, we conclude with a proof of the existence of solutions to the Euler-Lagrange equations.

*Proof of Theorem* 1.13. If we can verify Gâteaux differentiability of $F_s$, then the proof is completed, since defining

$$I(t) := F_s(u + t\varphi),$$

where $u$ is any minimizer of $F_s$ over $H_g^{s,p}(\Omega)$ and $\varphi \in C_c^\infty(\Omega)$, then $I$ is differentiable and

$$I(0) = \min\{I(t) : t \in \mathbb{R}\}.$$

Thus,

$$I'(0) = \frac{d}{dt} F_s(u + t\varphi) = \langle F_s'(u), \varphi \rangle.$$

It therefore remains to verify Gâteaux differentiability of $F_s$. However, we have

$$|\langle F_s'(u), \varphi \rangle| = \left| \int_{\mathbb{R}^d} f_z(x, u, D^s u) \, \varphi + D_\xi f(x, u, D^s u) \cdot D^s \varphi \, dx \right|$$

$$\leq \int_{\mathbb{R}^d} \left( C(|u|^{p-1} + |D^s u|^{p-1}) + \gamma_2(x) \right) |\varphi| \, dx$$

$$+ \int_{\mathbb{R}^d} \left( C(|u|^{p-1} + |D^s u|^{p-1}) + \gamma_3(x) \right) |D^s \varphi| \, dx$$

$$\leq \left\| C(|u|^{p-1} + |D^s u|^{p-1}) + \gamma_2(\cdot) \right\|_{L^{p'}(\mathbb{R}^d)} \|\varphi\|_{L^p(\mathbb{R}^d)}$$

$$+ \left\| C(|u|^{p-1} + |D^s u|^{p-1}) + \gamma_3(\cdot) \right\|_{L^{p'}(\mathbb{R}^d)} \|D^s \varphi\|_{L^p(\mathbb{R}^d)}$$

$$\leq \left( C(\|u\|_{L^p(\mathbb{R}^d)}^{\frac{p}{p'}} + \|D^s u\|_{L^p(\mathbb{R}^d)}^{\frac{p}{p'}}) + \|\gamma_2\|_{L^{p'}(\mathbb{R}^d)} \right) \|\varphi\|_{L^p(\mathbb{R}^d)}$$

$$+ \left( C(\|u\|_{L^p(\mathbb{R}^d)}^{\frac{p}{p'}} + \|D^s u\|_{L^p(\mathbb{R}^d)}^{\frac{p}{p'}}) + \|\gamma_3\|_{L^{p'}(\mathbb{R}^d)} \right) \|D^s \varphi\|_{L^p(\mathbb{R}^d)}$$

This show that $F$ is Gâteaux differentiable and the proof is complete.     □

**Acknowledgments.** The first author would like to thank Professor Ming-Chih Lai for his support during a portion of this project. The second author would like to thank Eliot Fried for discussions regarding non-local continuum mechanics and also for making him aware of the work of Edelen and Laws, Edelen, Green and Laws. Finally the authors would like to thank Armin Schikorra for many helpful discussions regarding the paper.

## Author information

Tien-Tsan Shieh, National Center for Theoretical Sciences, National Taiwan University, Taiwan.
E-mail: `ttshieh@ncts.ntu.edu.tw`

Daniel E. Spector, Department of Applied Mathematics, National Chiao Tung University, Taiwan.
E-mail: `dspector@math.nctu.edu.tw`